\documentclass[reqno]{amsart}
\usepackage[process=all,crop=pdfcrop]{pstool}
\usepackage{float}

\numberwithin{equation}{section} 
\usepackage{amsmath,amsfonts,amssymb,amsthm,latexsym}
\usepackage{enumerate}
\usepackage{cancel}
\usepackage{mathrsfs}
\usepackage{stackrel}

\usepackage{cases}
\usepackage[square,comma,numbers,sort&compress]{natbib}

\usepackage[colorlinks=true]{hyperref}
\hypersetup{urlcolor=blue, citecolor=red}

\newcounter{mnote}
\theoremstyle{plain}
\newtheorem{theorem}{Theorem}[section]

\theoremstyle{definition}

\theoremstyle{remark}

\newcommand{\field}[1]{\mathbb{#1}}
\newcommand{\nN}{\field{N}}

\newcommand{\nR}{\field{R}}

\newcommand{\OOO}{\mathcal O}

\newcommand{\sand}{\quad\text{and}\quad}

\newcommand{\pd}[2]{\frac{\partial #1}{\partial #2}}
\newcommand{\od}[2]{\frac{d #1}{d #2}}

\newcommand{\norm}[1]{\left\lVert#1\right\rVert}
\newcommand{\set}[1]{\left\{#1\right\}}

\newcommand{\LpP}[1]{\text{$L^{#1}$}}
\newcommand{\HpP}[1]{\text{$H^{#1}$}}

\renewcommand{\Pr}{\text{Pr}}
\newcommand{\Ra}{\text{Ra}}
\newcommand{\Nu}{\text{Nu}}
\newcommand{\tu}{\tilde {u}}

\newcommand{\vt}{\tilde {v}}
\newcommand{\ttheta}{\tilde {\theta}}

\newcommand{\tomega}{\tilde {\omega}}
\newcommand{\nc}{n_{\text C}}
\newcommand{\nf}{n_{\text F}}

\begin{document}
\title[Data assimilation of nearly turbulent Rayleigh-B\'enard flow]{Assimilation of nearly turbulent Rayleigh-B\'enard flow through vorticity or local circulation measurements: a computational study}

\date{\today}

%
\author{Aseel Farhat}
\address[A. Farhat]{Department of Mathematics\\
                University of Virginia,\\
       Charlottesville, VA 22904, USA}
\email {af7py@virginia.edu}
\author{Hans Johnston}
\address[H. Johnston]{Department of Mathematics and Statistics\\
                University of Massachusetts, \\
        Amherst, MA 01003, USA}
\email {johnston@math.umass.edu}
\author{Michael Jolly}
\address[M. S. Jolly]{Department of Mathematics\\
                Indiana University, Bloomington\\
        Bloomington , IN 47405, USA}
\email {msjolly@indiana.edu}
\author{Edriss S. Titi}
\address[E. S. Titi]{Department of Mathematics, Texas A\&M University, 3368 TAMU,
 College Station, TX 77843-3368, USA.  {Also},
 Department of Computer Science and Applied Mathematics, Weizmann Institute
 of Science, Rehovot 76100, Israel.} 
 \email{titi@math.tamu.edu and
  edriss.titi@weizmann.ac.il}

\begin{abstract}
We introduce a continuous (downscaling) data assimilation algorithm for the 2D B\'enard convection problem using vorticity or local circulation measurements only. 
In this algorithm, a nudging term is added to the vorticity equation to constrain the
model. Our numerical results indicate that the approximate solution of the algorithm is converging to the unknown reference solution (vorticity and temperature) corresponding to the measurements of the 2D B\'enard convection problem when only spatial coarse-grain measurements of vorticity are assimilated. Moreover, this convergence is realized using data which is much more coarse than the resolution needed to satisfy rigorous analytical estimates. 
\end{abstract}

 \maketitle
   \centerline{\it  This paper is dedicated to our great friend and collaborator Professor Bernardo Cockburn}
   \centerline{\it  on the occasion of his sixtieth birthday.}
   \vskip .15truein



\bigskip
\section{Introduction}\label{intro}

There are a variety of data assimilation algorithms whereby spatial coarse measurements are used to enhance the accuracy of model forecasts.
Ensemble Kalman filters incorporate statistical noise to estimate both the state and the uncertainty of the forecast \cite{Browning_H_K, B_L_Stuart, Evensen}. Though widely used, Kalman filters are mainly justified rigorously only for linear systems with Gaussian noise. Early work on continuous data assimilation by feedback nudging dates back to the 1970s. A direct method for data assimilation proposed by Charney, Halem and Jastrow \cite{chj} (see also \cite{B_L_Stuart, Olson_Titi_2003, Henshaw_Kreiss_Ystrom, Ghil_Halem_Atlas, Ghil_Shkoller_Yangarber}) continuously injects coarse grain observational data into
the nonlinear term in the model. Since this term typically involves spatial differentiation, implementation for most discretizations can be difficult.
Differentiation of the coarse-grain data can be avoided by inserting the data instead into a feedback nudging term.   
 
The algorithm used here was adapted in \cite{Azouani_Olson_Titi} from one for stabilizing solutions \cite{Azouani_Titi, Lunasin_Titi, CaoKevrekidisTiti}. 
It is easily described in terms of a general dissipative differential equation of the form
\begin{align}
&\od{u}{t} = F(u), \label{1.1}\\
&u(0) = u_0. \label{1.2}
\end{align}
If the initial data $u_0$ is given, then one can integrate \eqref{1.1}-\eqref{1.2} to find the corresponding solution. However, in weather forecasting and other applications, the initial data $u_0$ is unknown. The goal of data assimilation is to design an algorithm for recovering $u$ from coarse-grain spatial measurements over a certain interval of time. Specifically, for an unknown (reference) solution $u$ to \eqref{azouani}, for which in particular the initial data $u_0$ is unknown, 
we seek to recover that solution from spatial coarse-grain data through the associated downscaling data assimilation algorithm
\begin{subequations}\label{azouani}
\begin{align}
&\od{\tu}{t} = F(\tu) - \mu (I_h(\tu)- I_h(u)), \\
&\tu(0)= \tu_0,
\end{align}
\end{subequations}
where $\mu>0$ is a relaxation (nudging) parameter, $\tu_0$ is arbitrary,
and $I_h$ is some suitable approximation of the identity interpolant operator, constructed from the spatial coarse-grain measurement of $u$, with spatial resolution $h$.

The idea behind this algorithm, say in the case of the 2D Navier-Stokes equations (NSE), capitalizes on the fact that turbulent flows are determined by finitely many parameters, such as determining modes, nodes and local spatial averages \cite{C_O_T_1995, C_O_T_1997,Foias_Prodi, Foias_Temam_2, Foias_Titi, Jones_Titi}. In particular, the instabilities in such flows occur at the large spatial scales, while the fine scales are dissipated by the viscosity. Consequently, the algorithm suggested in \eqref{azouani} stabilizes, through the nudging term, the large spatial scales of size $h$ and above. While the nudging term is stabilizing the large scale, it might cause instabilities in the smaller spatial scales magnified by the nudging parameter $\mu>0$, a phenomenon known in control theory as ``spill over". This additional instability in smaller scales, however, can be dissipated by the viscosity; which is the main reason for limiting the size of the nudging parameter $\mu$, i.e., $\mu$ should be large, but not too large. The rigorous treatment of the algorithm \eqref{azouani} seeks to establish the right balance between $h$, to be small but not too small, and $\mu$ to be large, but not too large, to allow for the viscosity to stabilize the fine spatial scales. 

A variety of systems have been treated by this algorithm. For the incompressible 2D Navier-Stokes equations (NSE), it is shown in \cite{Azouani_Olson_Titi} that if $\mu$ is sufficiently large and $h$ sufficiently small, specified explicitly in terms of the physical parameters, then $\|v(t) - u(t) \|_{L^2} \to 0$ (or $\norm{v(t) - u(t)}_{H^1} \to 0$ under further assumptions on the size of $\mu$), at an exponential rate (see also the computational work in \cite{Altaf, Gesho_Olson_Titi}). The convergence of this synchronization algorithm for the 2D NSE, in higher order (Gevery class) norm and in $L^\infty$ norm, was later studied in \cite{B_M} for smoother forcing. An extension of the approach in \cite{Azouani_Olson_Titi} to the case when the observational data contains stochastic noise was analyzed in \cite{Bessaih_Olson_Titi}. A study of the algorithm for the 2D NSE when the measurements are obtained discretely in time and are contaminated by the measurement device error is presented in \cite{FMTi} (see also \cite{Hayden_Olson_Titi}).  Convergence of this nudging type algorithm has been proved for the 3D Brinkman-Forchheimer-Extended Darcy model  \cite{MTT2015}, the 3D Navier-Stokes-$\alpha$ model \cite{Albanez_Nussenzveig-Lopes_Titi}, the subcritical surface quasi-geostrophic (SQG) equation \cite{J_M_T} and the damped driven Kortevieg-de Vries (KdV) equation \cite{Jolly_Sadigov_Titi_2}. More recently in \cite{M_T}, the authors obtain uniform in time estimates for the error between the numerical approximation given by the Post-Processing Galerkin method of the downscaling algorithm and the reference solution, for the 2D NSE. This provides evidence and rigorous analytic justification that the numerical implementation of this algorithm is practically reliable. For certain models this approach can work with data in only a subset of the system state variables. This is proved for the 2D NSE in \cite{Farhat_Lunasin_Titi1}, and for the 3D Leray-$\alpha$ model of turbulence model in \cite{Farhat_Lunasin_Titi5} for data collected for only one and two component(s) of velocity, respectively. Notably, the treatment of each of the above mentioned systems has its own subtleties, and the studies are motivated by specific scientific questions, as we will clarify below. For instance, the treatment of the damped and driven KdV  \cite{Jolly_Sadigov_Titi_2}, capitalizes on the various conserved functionals of these systems in the absence of damping and forcing. This is needed because of the lack of regularization mechanism (no viscosity) in the KdV equation, unlike the other parabolic systems, such as the NSE.  

The Rayleigh-B\'enard system governs two fields, velocity and temperature.   It was proved in \cite{Farhat_Jolly_Titi},
that both velocity and temperature can be recovered using only velocity coarse-mesh spatial measurements. Later in \cite{Farhat_Lunasin_Titi3}, it was shown that a modification of the algorithm \eqref{DA_Bous} for the Rayleigh-B\'enard system, with no penetration stress-free boundary conditions on the horizontal walls, works by employing only coarse-mesh spatial measurements for the horizontal component of velocity vector field. 

The analytical results reported in \cite{Farhat_Jolly_Titi} were demonstrated numerically, for low Rayleigh numbers, in \cite{Altaf} where evidence 
was also provided to show that assimilation can fail when using temperature measurement data alone. The works in \cite{Farhat_Jolly_Titi, Altaf} were performed under no-slip boundary conditions for the velocity vector field at the horizontal solid walls. The question whether temperature observations are enough to determine all the dynamical state of the system is an important practical problem in meteorology and engineering and referred to as {\it Charney's conjecture} \cite{chj}. A precise theoretical formulation of the Charney's conjecture for certain simple atmospheric models is presented in \cite{Ghil_Shkoller_Yangarber}. The analytical argument in \cite{Ghil_Shkoller_Yangarber} suggested that the conjecture can be correct. On the other hand, the numerical tests in \cite{Ghil_Shkoller_Yangarber, Ghil_Halem_Atlas} for the primitive equations (and recently in \cite{Altaf} for the B\'enard convection problem) suggest that it is not certain that assimilation using coarse-grain temperature measurements alone will always recover the true state of the system. A recent study in \cite{Farhat_Lunasin_Titi2} shows that the data assimilation algorithm \eqref{azouani} using coarse-grain temperature measurement alone converges, at an exponential rate in time, to the reference solution of the B\'enard convection problem in porous medium. Similar results have also been established in \cite{Farhat_Lunasin_Titi4} for the 3D viscous Planetary Geostrophic circulation model. The results mentioned above provide rigorous justification for Charney's conjecture applied to certain simple planetary scales atmospheric models, while in practice Charney's conjecture remains open, and not certainly true, for other weather prediction models. 

In this paper we carry out a computational experiment similar to that in \cite{Altaf}, but at higher Rayleigh numbers and on the
vorticity formulation of the Rayleigh-B\'enard system. We propose, and computationally test, a data assimilation algorithm in the spirit of the algorithm in \cite{Azouani_Olson_Titi, Farhat_Jolly_Titi}, that uses coarse-grain vorticity or local circulation spatial measurements only rather than velocity measurements, to recover the full state of the system. The algorithm is presented in section \ref{algorithm} and the numerical results are presented in section \ref{numerical}. 

\bigskip

\section{Mathematical Analysis of the B\'enard Convection Problem}

The Rayleigh-B\'enard convection problem is a model of the Boussinesq system of an incompressible fluid layer, confined between two solid walls,  which is heated from below and cooled from the top in such a way that the lower wall maintains a temperature $T_0$ while the upper one maintains a temperature $T_1<T_0$. In this case, the two-dimensional Boussinesq equations that govern the velocity, pressure and temperature in $\Omega=[0,L]\times[0,1]$ can be written as (see \cite{Foias_Manley_Temam})
\begin{subequations}\label{Bous}
\begin{align}
&\pd{u}{t} - \nu\Delta u + (u\cdot\nabla)u + \nabla p = \theta \mathbf{e}_2, \label{Bous1}\\
&\pd{\theta}{t} - \kappa\Delta\theta + (u\cdot\nabla)\theta =0, \label{Bous2}\\
&\nabla\cdot u= 0,\label{Bous_div}\\
&u(0;x) = u_0(x), \quad \theta(0;x)=\theta_0(x), \label{Bous_initial}
\end{align}
where $\nu = (\Pr/\Ra)^{1/2}$ and $\kappa = (\Pr\Ra)^{-1/2}$, $\Pr$ is the Prandtl number and $\Ra$ is the Rayleigh number. 
We complement the above system with the following boundary conditions: 
\begin{align}
u = 0 \quad \text{at} \quad x_2=0 \quad \text{and} \quad u = 0 \quad \text{at} \quad x_2=1, \label{boundary1}\\
\theta = 1 \quad \text{at} \quad x_2=0 \quad \text{and} \quad  \theta = 0 \quad \text{at} \quad x_2=1, \label{boundary2}
\end{align}
and
\begin{align}
& u, \theta, p \text{ are periodic, of period } L, \text{ in the }x_1\text{-direction}.\label{boundary3}
\end{align}
\end{subequations}

The data assimilation algorithm proposed in \cite{Farhat_Jolly_Titi} for this problem is given by 
\begin{subequations}\label{DA_Bous}
\begin{align}
&\pd{\tu}{t} -\nu \Delta \tu + (\tu\cdot\nabla)\tu +\nabla \tilde p = \ttheta\mathbf{e}_2- \mu(I_h(\tu)-I_h(u)), \\
&\pd{\ttheta}{t} -\kappa\Delta\ttheta - (\tu\cdot\nabla)\ttheta= 0, \\
&\nabla \cdot \tu = 0, \\
&\tu(0;x) = \tu_0(x), \quad \ttheta(0; x) = \ttheta_0(x), \label{DA_Bous_initial}
\end{align}
\end{subequations}
with boundary conditions \eqref{boundary1} \eqref{boundary3} holding for 
$u, \theta$ and $p$  replaced by $\tu, \ttheta$ and a modified pressure $\tilde p$, respectively. Here $I_h$ is a linear interpolating operator at spatial resolution $h$. It is worth stressing that the algorithm \eqref{DA_Bous} employs only spatial coarse-grain measurements of the velocity field and it does not require any measurements of the temperature, an interesting feature of this method, as we will discuss below.

Two types of interpolating operators can be considered. One satisfies
\begin{align}\label{app}
I_h: \HpP{1} \rightarrow \LpP{2} \quad \text{and} \quad \norm{\xi - I_h(\xi)}_{\LpP{2}}^2 \leq c_0h^2\norm{\xi}_{\HpP{1}}^2, 
\end{align}
for every $\xi \in \HpP{1}$, where $c_0>0$ is a dimensionless constant.
The other satisfies
\begin{align}\label{app2}
I_h: \HpP{2}\rightarrow\LpP{2} \quad \text{and} \quad \norm{\xi - I_h(\xi)}_{\LpP{2}}^2 \leq c_0h^2\norm{\xi}_{\HpP{1}}^2 + c_0^2h^4\norm{\xi}_{\HpP{2}}^2,
\end{align}
for every $\xi \in \HpP{2}$, where $c_0>0$ is a dimensionless constant. 

A physical example of an interpolant observable that satisfies \eqref{app} is based on local volume elements (local volume averages) that were studied in \cite{Azouani_Olson_Titi, Jones_Titi}, see also discussion below. A physical example of an interpolant observable that satisfies \eqref{app2} is based on measurements at a discrete set of nodal points in $\Omega$ (see Appendix A in \cite{Azouani_Olson_Titi}).  Flexibility in the choice of interpolant is one of the main advantages of injecting the observed data thorough a feedback nudging term, as opposed to inserting it in the nonlinear term, which  in most models would require computing spatial derivatives of the coarse-grain measurements. 

It is shown in \cite{Farhat_Jolly_Titi} that \eqref{DA_Bous} has a unique solution for all $t \ge 0$.  Thus if we were to set $\tu_0=u_0$ and $\ttheta_0=\theta_0$, we
would have $(\tu,\ttheta)(t)=(u, \theta)(t)$ for all $t \ge 0$. Here, however,
we assume that exact initial data $u_0$ is unknown, and use spatially coarse-grain data for only the velocity vector field of a reference
solution to \eqref{Bous} to drive or {\it nudge} the solution of \eqref{DA_Bous} toward that reference solution
at an exponential rate.  No temperature data is used in the nudging.

Let  $\mathcal{F}$ be the set of $C^\infty(\Omega)$ functions defined in $\Omega$, which are trigonometric polynomials in $x_1$ with period $L$, and compactly supported in the $x_2$-direction. We denote the space of smooth vector-valued functions incorporating the divergence-free condition by
\[\mathcal{V}:=\set{\phi\in\mathcal{F}\times\mathcal{F}: \; \nabla\cdot\phi=0}.\]
The closures of $\mathcal{V}$ and $\mathcal{F}$ in $L^2(\Omega)$ will be denoted by $H_0$ and $H_1$, endowed with the usual scalar product
\[(u,v)_{H_0}=\sum_{i=1}^2\int_{\Omega} u_i(x)v_i(x)\,dx
\sand
(\psi,\phi)_{H_1}=\int_{\Omega} \psi(x)\phi(x)\,dx, \]
and the associated norms $\norm{u}_{H_0} = (u,u)_{H_0}^{1/2}$ and $\norm{\phi}_{H_1} = (\phi,\phi)_{H_1}^{1/2}$, respectively.  Likewise we denote by  $V_0$ and $V_1$, the closures of $\mathcal{V}$ and $\mathcal{F}$ in $H^1(\Omega)$ respectively. The spaces $V_0$ and $V_1$ shall be endowed with the scalar product
\[((u,v))_{V_0}=\sum_{i,j=1}^2\int_{\Omega}\partial_ju_i(x)\partial_jv_i(x)\,dx
\sand
((\psi,\phi))_{V_1}=\sum_{j=1}^2\int_{\Omega}\partial_j\psi(x)\partial_j\phi(x)\,dx, \]
with associated norms $\norm{u}_{V_0} = ((u,u))_{V_0}^{1/2}$ and $\norm{\phi}_{V_1} = ((\phi,\phi))_{V_1}^{1/2}$, respectively.

It was shown in \cite{Foias_Manley_Temam, Temam_1997} that  the 2D B\'enard convection system has a finite-dimensional global attractor.

\begin{theorem}[Existence of a global attractor]\cite{Foias_Manley_Temam, Temam_1997}\label{global_attractor_Bous} Let $T>0$ be fixed, but arbitrary. If the initial data $u_0\in V_0$ and $\theta_0\in{V_1}$, then the system \eqref{Bous} has a unique strong solution $(u,\theta)$ that satisfies $u\in C([0,T];V_0)\cap L^2([0,T];\mathcal{D}(A_0))$ and $\theta\in C([0,T];{V_1})\cap L^2([0,T];\mathcal{D}(A_1))$.  Moreover the system induced by \eqref{Bous} is well-posed and possesses a finite-dimensional global attractor $\mathcal{A}$ which is maximal among all the bounded invariant, for all $t\in\nR$, sets and is compact in $H_0\times {H_1}$ and bounded in $V_0 \times V_1$.
\end{theorem}

The $V_0 \times V_1$ bounds on $(u,\theta)$ in the global attractor derived in \cite{Foias_Manley_Temam, Temam_1997} are exponential:
\begin{align}\label{Jbounds}
\|u\|_{V_0}^2, \|\theta \|_{V_1}^2 \le \rho= ae^b \quad  \;,
\end{align}
where $a=\OOO(\nu^{-3}),\ b=\OOO(\nu^{-8})$
in the case where the Prandtl number is near unity, i.e.,  $\Pr=\nu/\kappa \sim 1$.

\subsection{Rigorous Convergence Results}

The solution to \eqref{DA_Bous} converges to the reference solution as $t \to \infty$ at an exponential rate, provided
$\mu$ is sufficiently large, and accordingly, $h$ is sufficiently small.  The condition on $\mu$ and in turn on $h$ involves upper bounds on the global attractor of \eqref{Bous}.   

\begin{theorem}\cite{Farhat_Jolly_Titi} \label{th_conv_1}
Let $I_h$ satisfy the approximation property \eqref{app} and $(u(t),\theta(t))$ be a strong solution in the global attractor of 
\eqref{Bous}. If 
\begin{align}\label{mu_1}
\mu\gtrsim \frac{1}{\kappa\lambda_1} + \frac{\rho}{\nu}+ \frac{\rho^2}{\kappa^2\lambda_1\nu},
\end{align}
 and $h>0$ satisfies  
 \begin{align}\label{hcond}
 \mu c_0^2h^2 \le \nu\;,
 \end{align}
  where $\lambda_1=2\pi\min\{1,L^{-1}\}$,  then 
 $$
\norm{u(t)-\tilde{u}(t)}_{{H_0}}^2 + \norm{\theta(t)-\tilde{\theta}(t)}_{{H_1}}^2\le \exp(-\lambda_1 \min\{\nu,\kappa\}t)\;.
$$ 
\end{theorem}
A similar result is proved in \cite{Farhat_Jolly_Titi} for interpolating operators of the type in \eqref{app2}.

We present here an outline of a simplified proof, based on that in \cite{Farhat_Jolly_Titi}, but using the classic Gronwall lemma.
The key is to extract some dissipation from the feedback term by applying \eqref{app} along with \eqref{hcond}
to obtain 
\begin{equation}\label{extract}
\begin{aligned}
-\mu(I_h(w),w) = -\mu(I_h(w)-w, w) - \mu \norm{w}_{{H_0}}^2
&\leq \frac{\mu c_0^2h^2}{2}\norm{w}_{{V_0}}^2 - \frac{\mu}{2}\norm{w}_{{H_0}}^2 \\
 &\leq \frac{\nu}{2}\norm{w}_{{V_0}}^2 -\frac{\mu}{2}\norm{w}_{{H_0}}^2 \;.
\end{aligned}
\end{equation}
Then, setting $w=u-\tilde{u}$, $\xi=\vt-\tilde{\vt}$, we take the $H_0$ and $H_1$ scalar products of the equations for $w$ and $\xi$ respectively.  Applying
\eqref{extract}, we ultimately obtain
\begin{align*}
 \od{}{t} \left(\norm{w}_{{H_0}}^2 + \norm{\xi}_{{H_1}}^2\right) + \lambda_1\min\{\nu,\kappa\}\left(\norm{w}_{H_0}^2 + \norm{\xi}_{{H_1}}^2\right) \leq \\
\left[\frac{4}{\kappa\lambda_1} + \frac{c}{\nu}\norm{u}_{{V_0}}^2 + \frac{c}{\nu\kappa^2\lambda_1}\norm{\theta}_{{V_1}}^4 - \mu\right]\norm{w}_{{H_0}}^2\;,
\end{align*}
where $c$ is a universal, dimensionless constant.
Thus, if $\mu$ is large enough to make the bracketed expression $\le 0$, we have the desired exponential decay.

The resolution needed is then determined by the two conditions
$$
\frac{1}{\kappa\lambda_1} + \frac{1}{\nu}\norm{u}_{{V_0}}^2 +
\frac{1}{\nu\kappa^2\lambda_1}\norm{\theta}_{{V_1}}^4 \lesssim \mu \quad \text{and} \quad \mu h^2\lesssim \nu \;.
$$
For $\Pr\sim1$,  this would mean that in terms of the Rayleigh number $\Ra=(\nu\kappa)^{-1}$ we would need {\it at least} 
$$
\mu \gtrsim  \Ra^{3/2} \norm{\theta}_{{V_1}}^4\;,
$$
which combined with the bounds in \eqref{Jbounds}  requires that
$$
\mu \gtrsim  Ra^{9/2} e^{Ra^{4}} \qquad \text {and hence} \qquad   h \lesssim Ra^{-9/4} e^{-Ra^{4}}  \;.
$$
Note that even if the bounds in \eqref{Jbounds} could be reduced to $\OOO(1)$,
the last term in \eqref{mu_1} would still require that $\mu \gtrsim \Ra^{3/2}$.   In the remaining sections we present numerical computations which demonstrate that this approach to data assimilation can be effective using much coarser data, and drastically smaller values for $\mu$, than those suggested by the analytical estimates above. 

\bigskip 

\section{A Data Assimilation Algorithm using Vorticity or local circulation spatial Coarse Mesh Measurements}\label{algorithm}
In vorticity-streamfunction formulation, the Rayleigh-B\'enard convection problem \eqref{Bous} can be written 
\begin{subequations}\label{Bous_Vor}
\begin{align}
&\pd{\omega}{t} - \nu\Delta \omega + (u\cdot\nabla)\omega = -\pd{\theta}{x_1}, \label{Bous1_Vor}\\
&\pd{\theta}{t} - \kappa\Delta\theta + (u\cdot\nabla)\theta =0, \label{Bous2_Vor}\\
&\omega(0;x) = \omega_0(x), \quad \theta(0;x)=\theta_0(x).\label{Bous_initial_Vor}
\end{align}
The no-slip boundary condition for the velocity field in \eqref{boundary1} and the divergence free condition \eqref{Bous_div} are enforced by the 
solution of an elliptic system for the streamfunction $\psi$,
given by
\begin{align} \label{BoussPsi}
\Delta \psi=\omega, 
\end{align}
with boundary conditions 
\begin{align}
\quad \psi =\frac{\partial \psi}{\partial x_2}=0\quad \text{at} \quad x_2=0 \quad \text{and} \quad x_2=1, \label{boundary1_Vor}
\end{align}
and velocity field defined in terms of $\psi$ by
$$u = \nabla^\perp \psi = (-\partial \psi / \partial x_2,\partial \psi / \partial x_1)^T. $$
The boundary condition \eqref{boundary3} is enforced by 
\begin{align}
& \omega, \psi, \theta, \text{ are periodic, of period } L, \text{ in the }x_1\text{-direction}, \label{boundary3_Vor}
\end{align}
with $\theta$ satisfying the boundary condition
\begin{align}
\theta = 1 \quad \text{at} \quad x_2=0 \quad \text{and} \quad  \theta = 0 \quad \text{at} \quad x_2=1. \label{boundary2_Vor}
\end{align}
\end{subequations}

The data assimilation algorithm that we propose and study here is given by: 
\begin{subequations}\label{Bous_Vor_DA}
\begin{align}
&\pd{\tomega}{t} - \nu\Delta \tomega + (\tu\cdot\nabla)\tomega = -\pd{\ttheta}{x_1} - \mu \left(I_h(\tomega) - I_h(\omega) \right), \label{Bous1_Vor_DA}\\
&\pd{\ttheta}{t} - \kappa\Delta\ttheta + (\tu\cdot\nabla)\ttheta=0, \label{Bous2_Vor_DA}\\
&\tomega(0;x) = \tomega_0(x), \quad \ttheta(0;x)=\ttheta_0(x),\label{Bous_initial_Vor_DA}
\end{align}
where 
\begin{align} 
\tomega = \Delta \phi, \quad \tu = \nabla^\perp \phi = (-\partial \phi / \partial x_2,\partial \phi / \partial x_1)^T,  
\end{align}
\end{subequations}
subject to the boundary conditions \eqref{boundary1_Vor}, \eqref{boundary3_Vor} and \eqref{boundary2_Vor} with $\psi$, $\omega$, and $\theta$ replaced by $\phi$, $\tomega$, and $\ttheta$, respectively. Here, $(\tomega_0,\ttheta_0)$ can be taken arbitrarily and $I_h$ is a linear interpolant operator constructed from the error-free discrete spatial measurements of the vorticity $\omega$.

Notice that when solving the B\'enard convection problem in vorticity-streamfunction formulation \eqref{Bous_Vor} (and similarly the data assimilation system \eqref{Bous_Vor_DA}), the streamfunction $\psi$ is subject to two boundary conditions, the Dirichlet and Neumann boundary conditions \eqref{boundary1_Vor}, making the elliptic equation \eqref{BoussPsi} overdetermined. On the other hand, the full system \eqref{Bous_Vor} (and similarly system \eqref{Bous_Vor_DA}) has the right number of required boundary conditions but it is not supplemented with explicit boundary conditions for the vorticity $\omega$, making equation \eqref{Bous1_Vor} (and similarly equation \eqref{Bous1_Vor_DA}) underdetermined. This is the same issue we have with the boundary conditions for the 2D Navier-Stokes equation in vorticity-streamfunction formulation. However, it is important to assign correct boundary conditions for the vorticity in order to obtain accurate assimilations of the flow. This issue with the boundary conditions has been treated in \cite{Ben-Artzi} for the 2D Navier-Stokes equations where ``the vorticity projection method" was introduced. The scheme was analyzed (a stability condition was derived) and tested numerically for several examples. The computational results in \cite{Ben-Artzi} indicate that ``the vorticity projection method" can be used to avoid the difficulty of determining the vorticity on the boundary. This method can be modified to be used for the 2D B\'enard convection system \eqref{Bous_Vor} and for the data assimilation algorithm \eqref{Bous_Vor_DA}.  

For the B\'enard system \eqref{Bous_Vor}, examples of interpolant observables that satisfy \eqref{app} can be taken as follows. One example is the projection on the low-modes of the Laplacian operator subject to periodic boundary conditions in the $x_1$-direction and the boundary conditions \eqref{boundary1_Vor} at the solid walls. That is, the projection $\Pi_h$ on the orthogonal eigenfunctions $\varphi_k$, with wave numbers $k\in \nN^2$ with $|k|\leq 1/h$, corresponding to the eigenvalues $\lambda_k \in \nR^{-}$, of the problem: 
\begin{align*}
\Delta \varphi_k &= \lambda_k \varphi_k, \\
\varphi_k(0, x_2) &= \varphi_k(L,x_2), \quad x_2 \in [0,1], \\
\varphi_k(x_1, 0) &= \varphi_k(x_1,1) = 0,  \quad x_1 \in [0,L].
\end{align*}
Notice that eigenfunctions $\{ \varphi_k(x_1,x_2): \, k\in \nN^2\}$ span the streamfunction $\psi(x_1,x_2)$ and $\norm{\Pi_h \psi - \psi}_{\LpP{2}} + \norm{\Pi_h \psi - \psi}_{\HpP{1}} \rightarrow 0,$ as $h\rightarrow 0$. The same eigenfunctions $\{ \varphi_k(x_1,x_2): \, k\in \nN^2\}$ also span the vorticity $\omega = \Delta \psi$ and $\norm{\Pi_h \omega - \omega}_\LpP{2}$, as $h \rightarrow 0$, though, the expansion is not expected to converge in $\HpP{1}$ since $\omega$ cannot be assumed to satisfy the Dirichlet boundary conditions in the $x_2$-direction. 

 Observe that one has to slightly modify the presentation given in \cite{Azouani_Olson_Titi, Jones_Titi} to fulfill the spatial boundary conditions \eqref{boundary1_Vor} -- \eqref{boundary3_Vor}. In the case of local volume elements of vorticity, we divide the domain $[0,L]\times[0,1]$ into $\{Q_j\}_{j=1}^N$ squares and the measurements of vorticity in each square are given by
$$ \overline{\omega_j} =  \frac{1}{|Q_j|}\int_{Q_j} \omega \, dx_1dx_2= \frac{1}{|Q_j|}\int_{Q_j} \nabla\times u \, dx_1dx_2.$$
By Green's Theorem in two-dimensions, 
$$ \overline{\omega_j} = \frac{1}{|Q_j|}\oint_{\partial Q_j} u \cdot \, d{\bf x} = \frac{1}{|Q_j|}\oint_{\partial Q_j} u_1dx_1 + u_2dx_2,$$
which is the local circulation at the boundary of $Q_j$. Notice that the local circulation is a quantity which might be measured while usually it is not easy to measure the vorticity. In that case, our algorithm can be implemented using local circulation measurements instead.

In the computations presented in this paper, we take $I_h = P_h$ to be the orthogonal projection onto the low Fourier modes in the $x_1$-direction and low Chebyshev modes in the $x_2$-direction, with wave numbers $k$ such that $|k|\leq 1/h$. 
Similar to the case of projection $\Pi_h$ discussed above, the expansion in the Fourier $\times$ Chebyshev modes, in the $x_1\times x_2$ directions, spans both the streamfunction and the vorticity functions $\psi(x_1,x_2)$ and $\omega(x_1,x_2)$ ($\phi$ and $\tomega$, respectively). On the boundary, the expansion in the Fourier $\times$ Chebyshev modes for the streamfunction $\psi$ (and $\phi$) will converge to the boundary conditions \eqref{boundary1_Vor} ($\psi, \phi = 0$ at $x_2 = 0,1$), but the expansion for the vorticity $\omega$ (and $\tomega$) might not converge at the boundary. Thus, $\norm{P_h \psi - \psi}_\LpP{2} + \norm{P_h \psi - \psi}_{\HpP{1}} \rightarrow 0$, as $h\rightarrow 0$, while $\norm{P_h \omega - \omega}_{\LpP{2}} \rightarrow 0$, as $h\rightarrow 0$, but the convergence is not expected to hold in $\HpP{1}$. In our computations, the values of $P_h(\omega)$ (and $P_h(\tomega)$) on the boundaries $x_2=0$ and $x_2=1$ correspond to the velocity field being no-slip, i.e., the streamfunction satisfies the Dirichlet boundary conditions in the $x_2$-direction.

An extension of the algorithm \eqref{Bous_Vor_DA} for the case of measurements with stochastic noise and the case of discrete spatio-temporal measurements with error can be established following the work in \cite{Bessaih_Olson_Titi} and \cite{FMTi}, respectively. 

A convergence theorem for the proposed algorithm \eqref{Bous_Vor_DA} can be established as in \cite{Farhat_Jolly_Titi}. The proof follows a similar argument to that for Theorem 3.3 in \cite{Farhat_Jolly_Titi}.  It is, however, more technical in this case so we do not present it in this paper in order to focus on the main goal, i.e., testing the practicality of the data assimilation algorithm \eqref{Bous_Vor_DA} numerically.  The analytical study of algorithm \eqref{Bous_Vor_DA} is a subject of future work. 

We have shown how the analytical lower bound on $\mu$ depends on the estimates on the global attractor from  \cite{Foias_Manley_Temam, Temam_1997}  which are exponential in the Rayleigh number $\Ra$. The numerical simulations for the 2D NSE in \cite{Gesho_Olson_Titi} (see also \cite{Hayden_Olson_Titi}) have shown that the algorithm \eqref{azouani} is effective in capturing the reference solution using much coarser data
than analytical estimates in \cite{Azouani_Olson_Titi} would suggest. For the 2D 
Rayleigh-B\'enard problem, it was shown in \cite{Altaf} that the reference solution can be recovered using much coarser velocity-only data
than required by the rigorous analysis in \cite{Farhat_Jolly_Titi}. Determining just how coarse the vorticity data in the algorithm \eqref{Bous_Vor_DA} can be in order to recover the reference solution of \eqref{Bous_Vor} motivates the computational study in this paper.

\bigskip

\section{Computational Results}\label{numerical}
%

\begin{figure}[h!] 
\centerline{\includegraphics[width=8.0cm]
{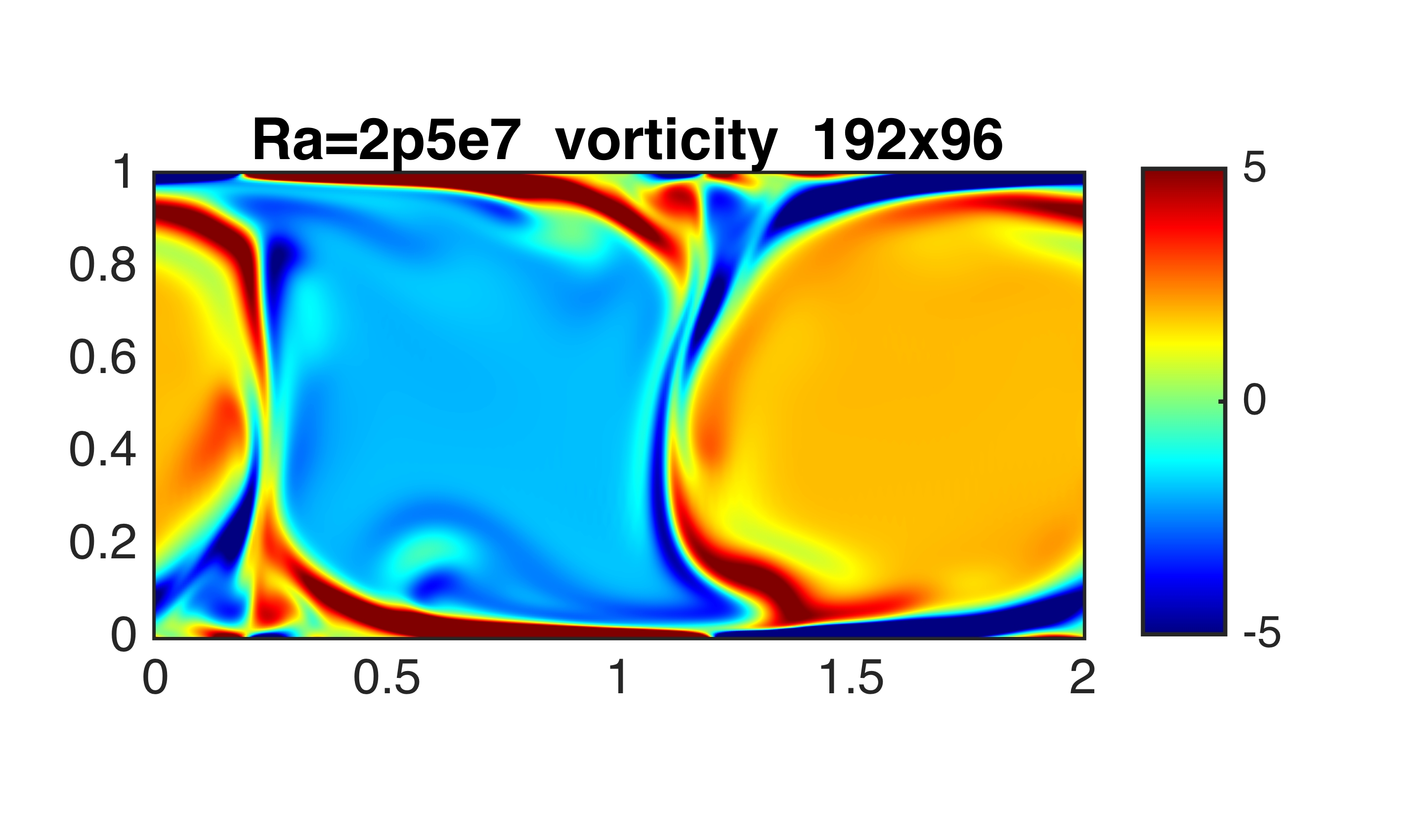} \hskip -14pt
\includegraphics[width=8.0cm]
{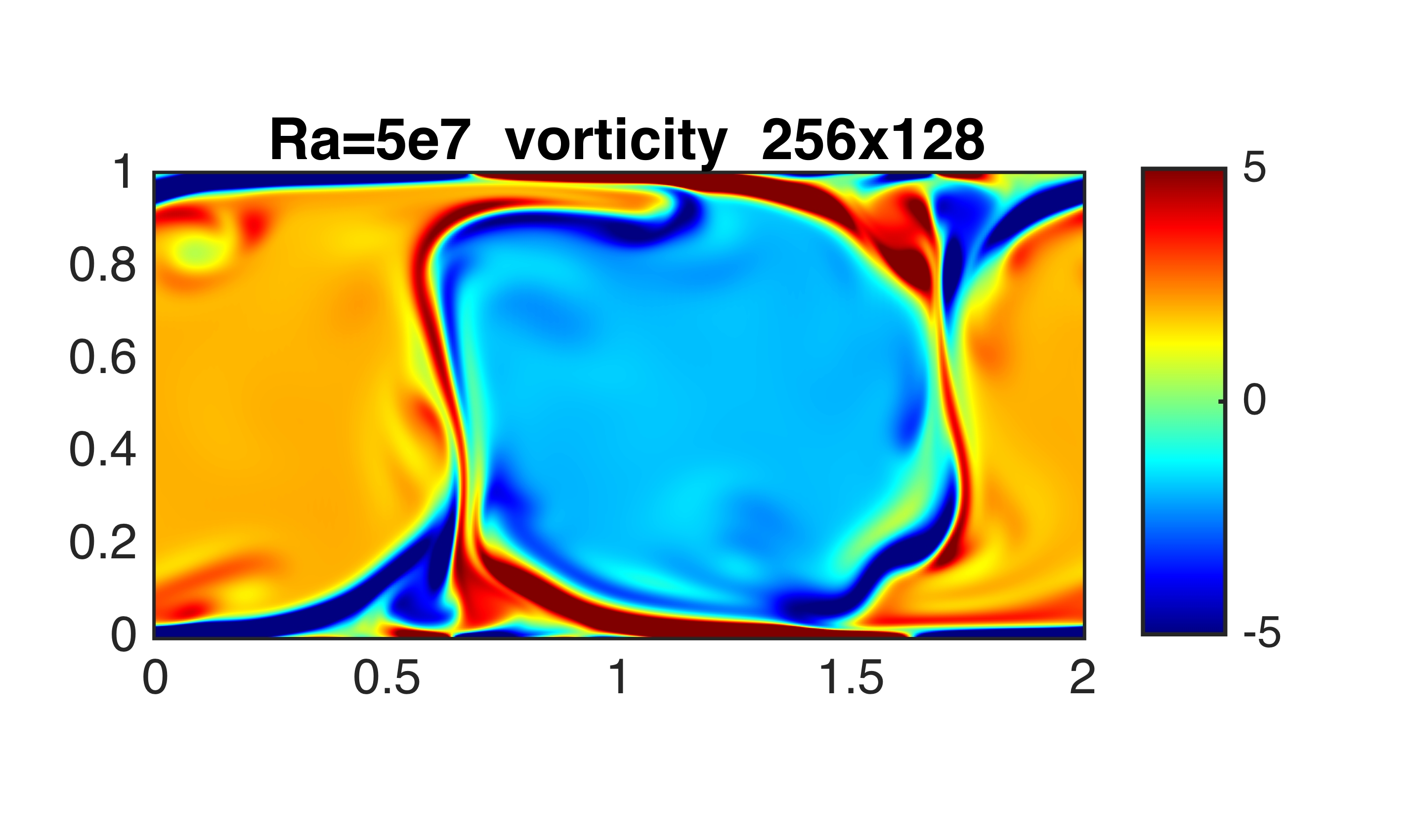}}
\vskip -20pt
\centerline{\includegraphics[width=8.0cm]
{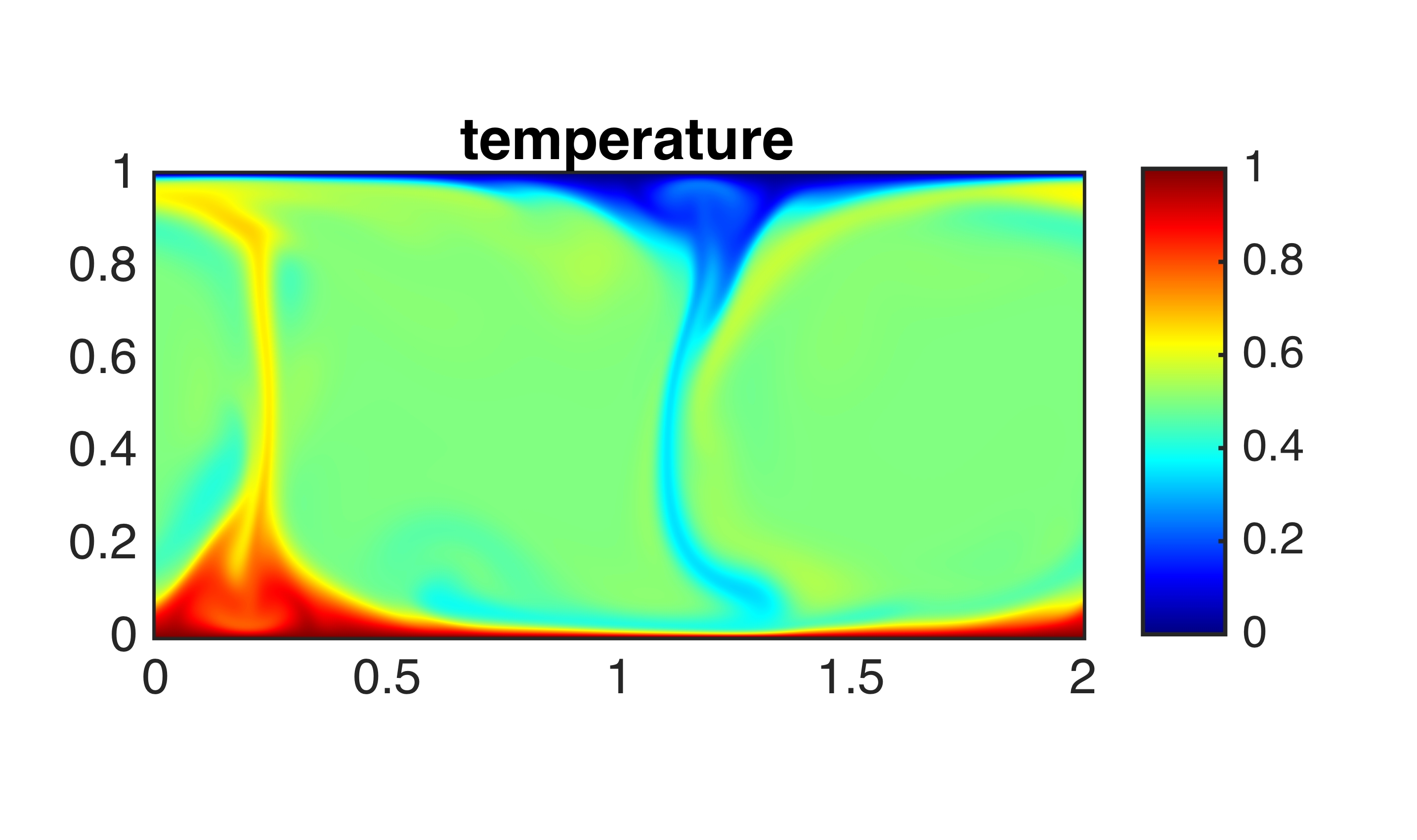} \hskip -14pt
\includegraphics[width=8.0cm]
{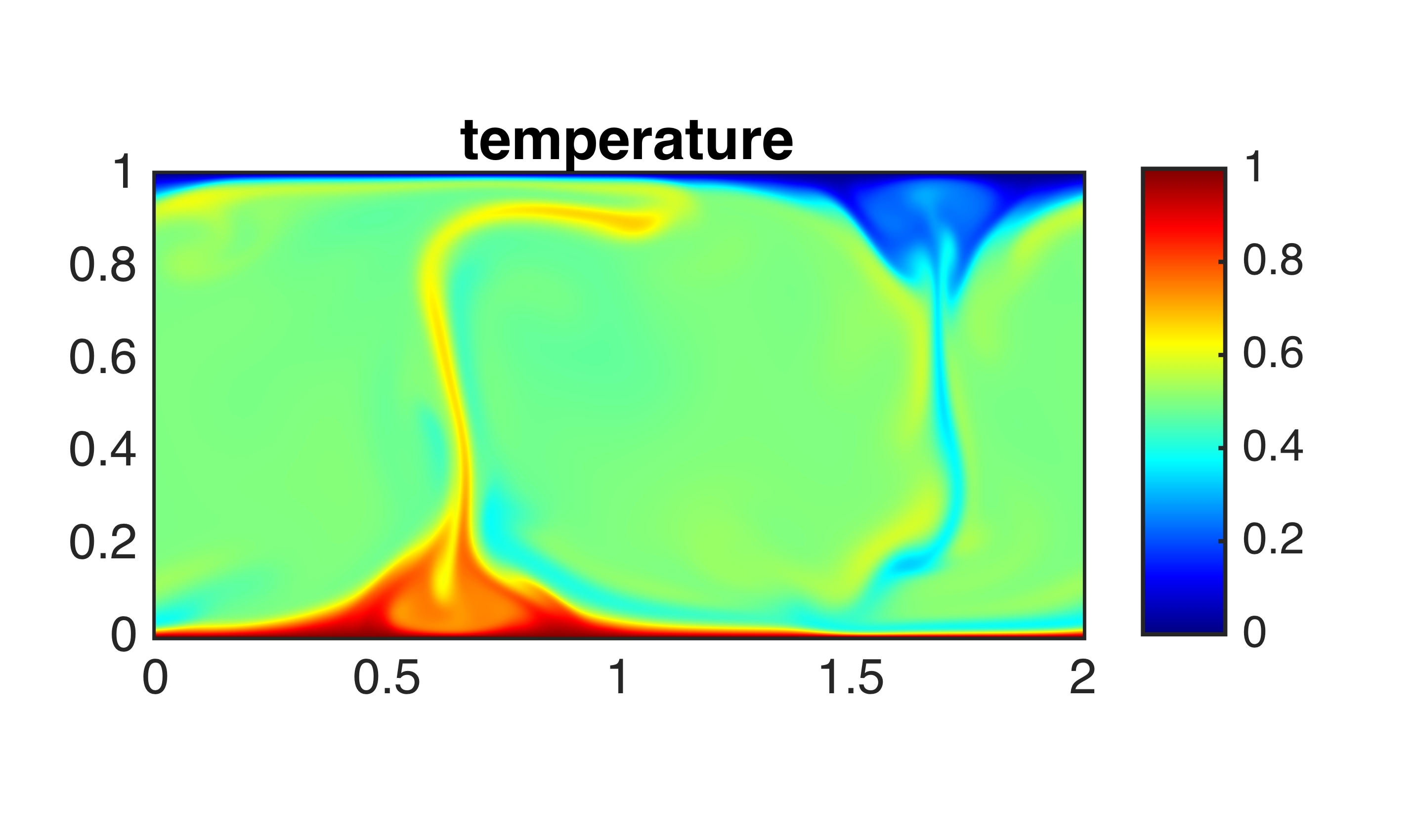}}
\caption{Initial data for reference solution}
\label{figure1}
\end{figure}

\begin{figure}[h!] 
\centerline{\includegraphics[width=7.0cm]
{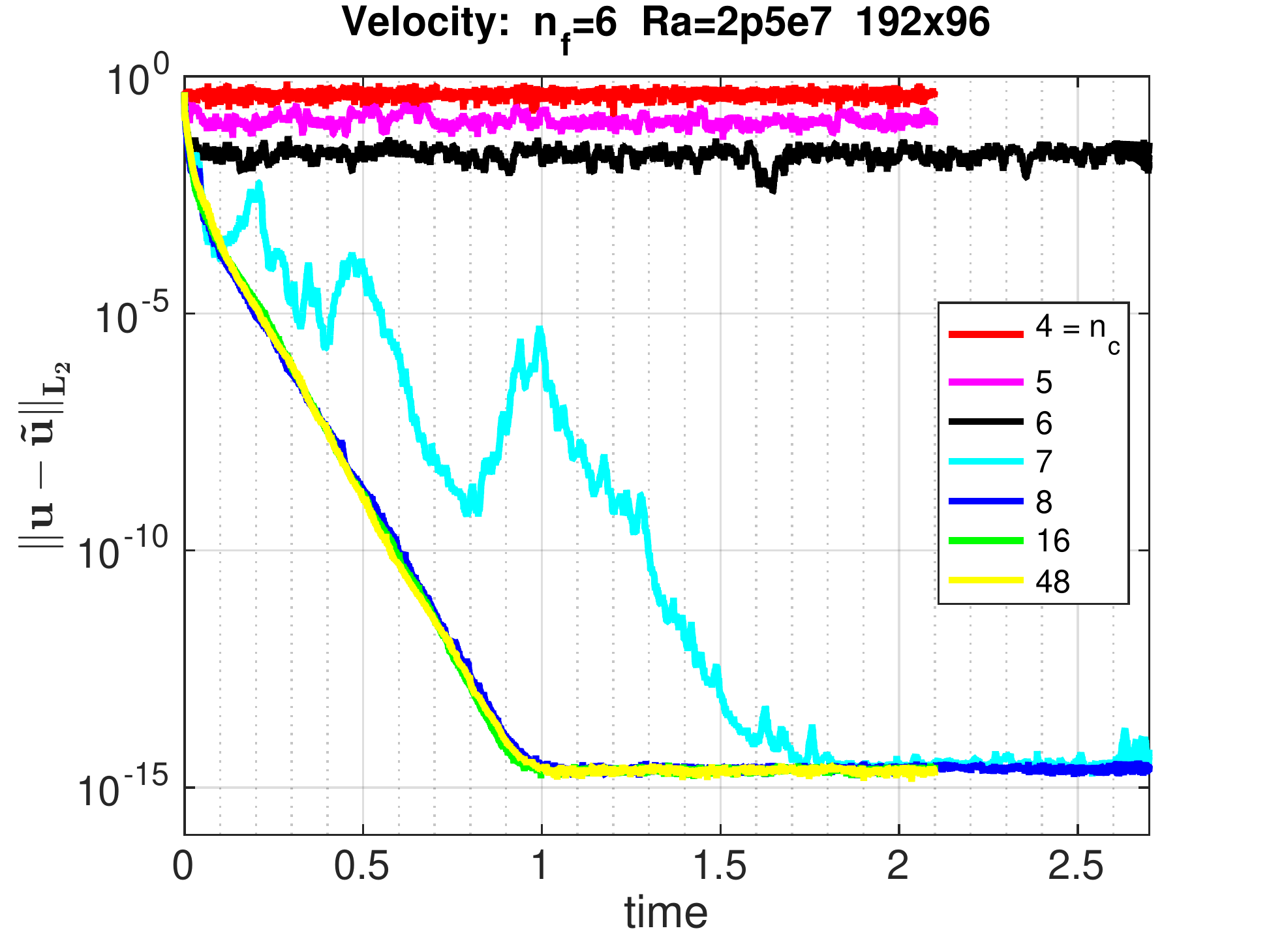} \
\includegraphics[width=7.0cm]
{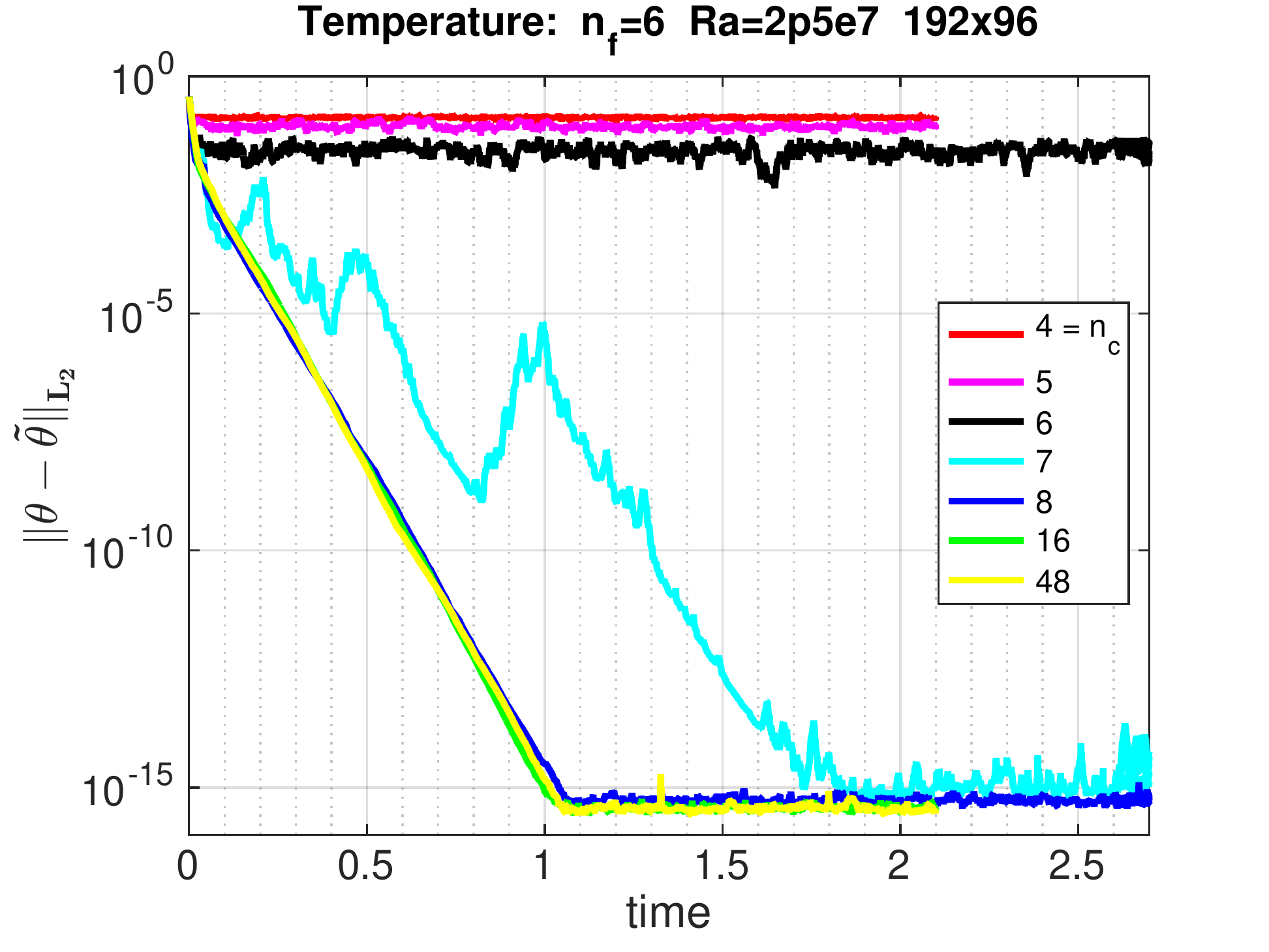}}
\centerline{\includegraphics[width=7.0cm]
{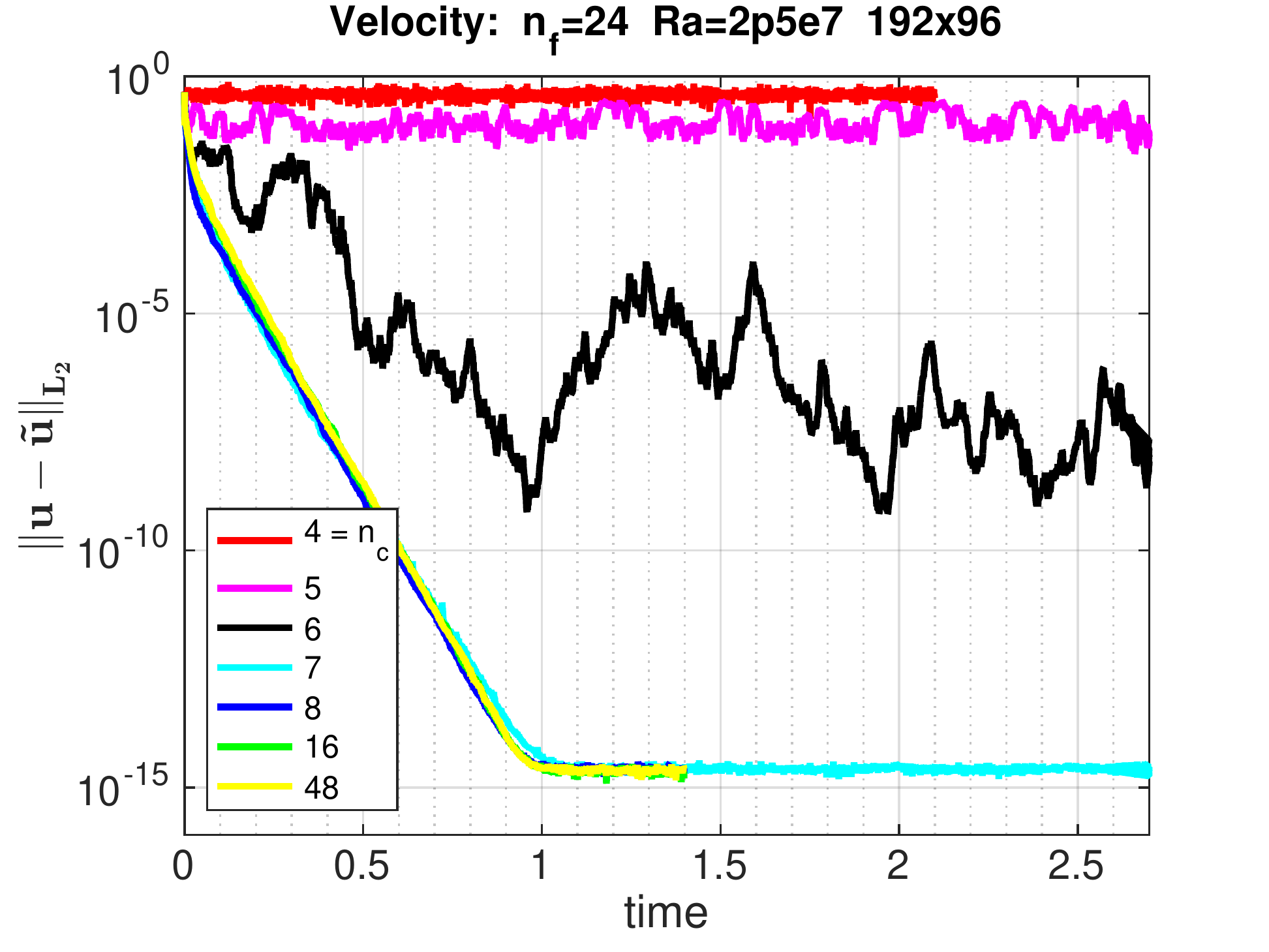} \ 
\includegraphics[width=7.0cm]
{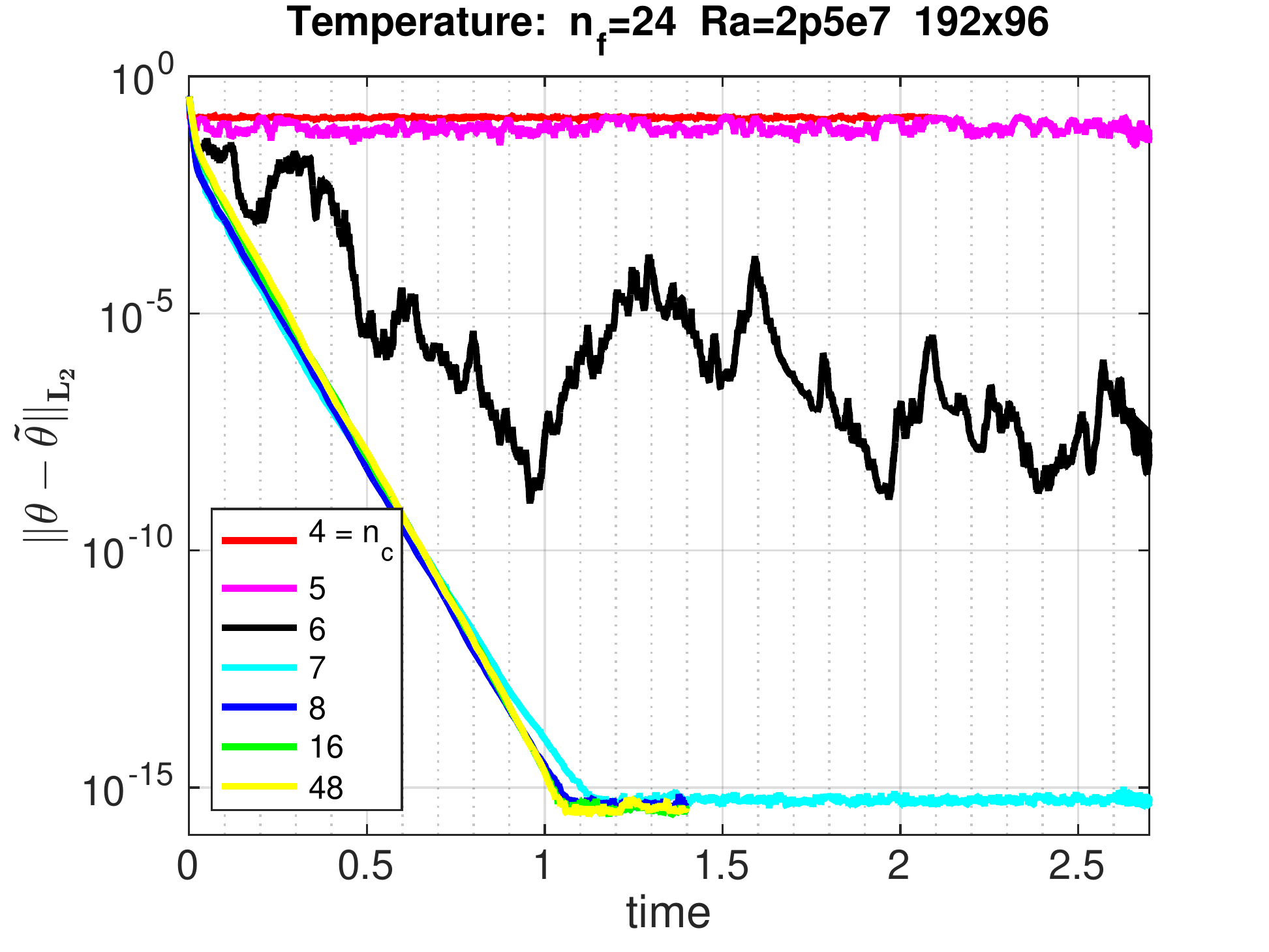}}
\caption{Convergence to reference solution at Ra$=2.5\times 10^7$}
\label{figure2}
\end{figure} 

We demonstrate that both the vorticity and temperature of the reference solution $(\omega,\theta)$ of \eqref{Bous_Vor}
can be recovered at an exponential rate from the solution $(\tomega, \ttheta)$ of the data assimilation algorithm \eqref{Bous_Vor_DA}. 
%
Here the role of the interpolant $I_h$ is played by $P_h$, the orthogonal projection onto the low Fourier modes in the $x_1$-direction and low Chebyshev modes in the $x_2$-direction, with wave numbers $k$ such that $|k|\leq 1/h$. We take the initial data
\begin{align}
\tu(0;x) = \tomega(0;x) =0, \quad \ttheta(0;x)=1-x_2, \ \forall  \ x \in [0,L] \times [0,1], 
\end{align}
with $L=2$. 

The numerical scheme for simulating \eqref{Bous_Vor} and \eqref{Bous_Vor_DA} is a 
Fourier-Chebyshev spectral collocation method in space and classical fourth order Runge-Kutta  
for the time stepping. Computation of the momentum equation and the elliptic 
equation \eqref{BoussPsi} are decoupled using a high order local formula for the 
vorticity at the boundary, derived from the Neumann boundary condition 
for the stream function in \eqref{boundary1_Vor}. The Dirichlet boundary condition
is enforced by the solution of the elliptic system 
 in \eqref{boundary1_Vor} which is solved by the matrix-diagonalization procedure.

\begin{figure}[h!] 
\centerline{\includegraphics[width=7.0cm]
{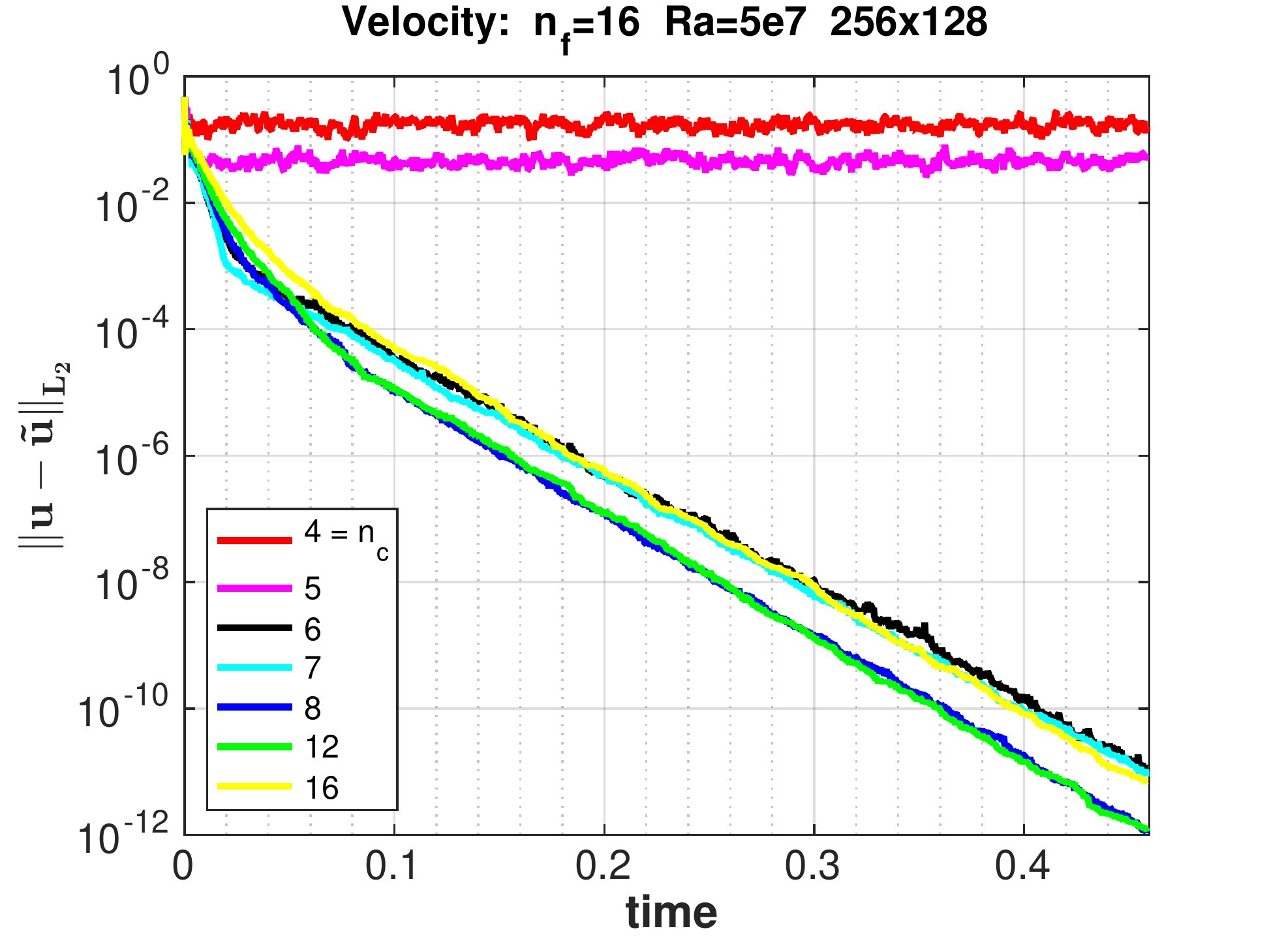} \
\includegraphics[width=7.0cm]
{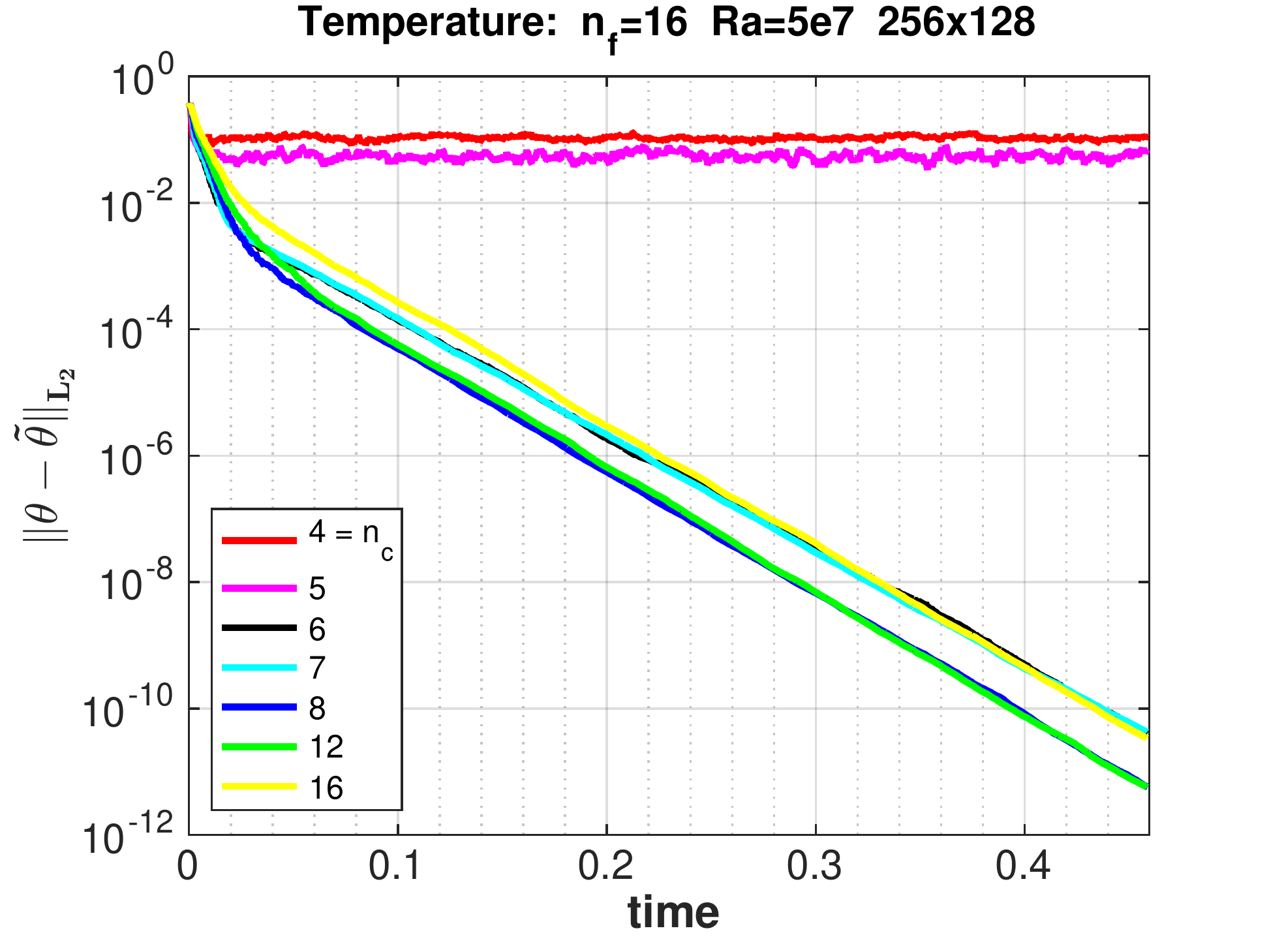}}
\centerline{\includegraphics[width=7.0cm]
{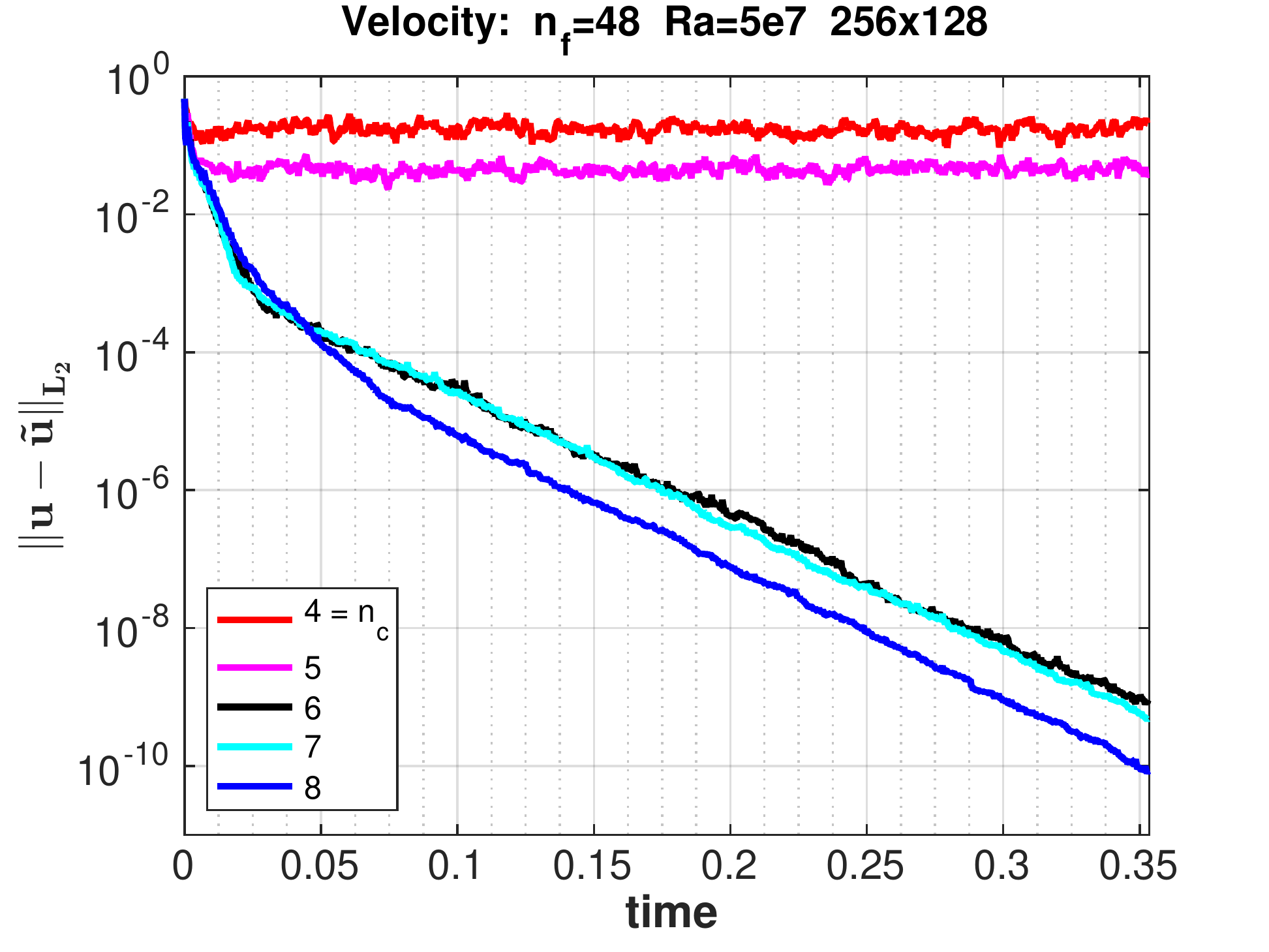} \ 
\includegraphics[width=7.0cm]
{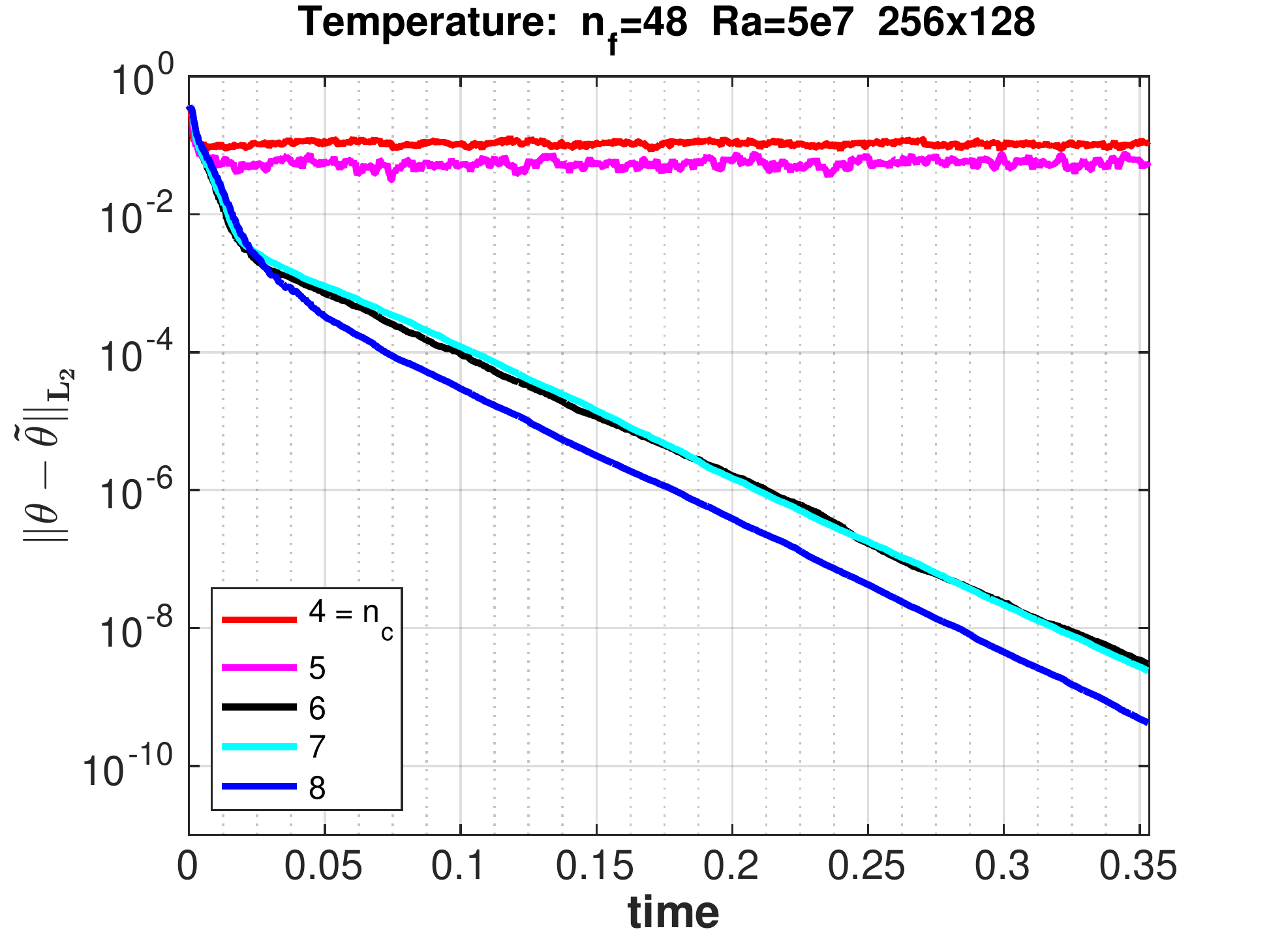}}
\caption{Convergence to reference solution at Ra$=5\times 10^7$}
\label{figure3}
\end{figure}

We present convergence results at two values of the Rayleigh number, 
$2.5\times 10^7$ and $5\times 10^7$, for which the reference solution 
would be considered to be in a moderately turbulent regime.  In the context
of the Nusselt number, defined by 
\begin{equation} \label{Nusselt}
\Nu = 1+{(\Pr\Ra)}^{1/2} {\langle (u\cdot{\bf e}_2) \theta\rangle },
\end{equation}
where $\langle \cdot \rangle$ denotes the space-time average,
this regime corresponds to $\Nu \approx 15$. 
$\Nu$ is a dimensionless quantity measuring the enhancement 
of the vertical heat transport by the convectively driven flow as
compared to the pure conduction solution.
Both the relaxation parameter and the Prandtl number are fixed at $\mu = \Pr = 1$. A snapshot of the initial condition for 
the reference solution for each Rayleigh number is shown in Figure \ref{figure1}.  
The time evolution of the $L^2$-norms for the difference in velocity and temperature 
for solutions to \eqref{Bous_Vor} and \eqref{Bous_Vor_DA} are shown in Figures \ref{figure2} and \ref{figure3}. The runs are
made for projectors $P_h$ with various resolutions to determine the minimum needed.  

At $\Ra=2.5\times 10^7$ we use a projection $P_h$ onto 6 and 24 
Fourier modes ($\nf$) in the $x_1$-direction, and from 4 to 48 Chebyshev modes ($\nc$) in the $x_2$-direction.
The solutions to both \eqref{Bous_Vor} and \eqref{Bous_Vor_DA} are computed by a pseudo-spectral code with $\nf\times\nc=192\times 96$ modes in which the Fourier modes are dealiased by a factor of 3/2.
Taking just $\nf\times \nc=6\times 8$ in $P_h$ suffices to nudge
the solution to \eqref{Bous_Vor_DA} to the reference solution at an exponential rate.  With $\nc=7$ in $P_h$ the solutions ultimately converge at an exponential rate, after a meandering 
period, while there is no convergence with fewer than 7 Chebyshev modes.
Using $\nf=24$ in $P_h$ enables an exponential rate convergence to begin immediately with $\nc=7$, but the borderline case of 6 Chebyshev modes 
meanders over the length of the run.

Meandering borderline cases are not seen at $\Ra=5 \times 10^7$, where we
use $\nf\times\nc=256\times 128$ dealiased modes to compute the reference solution.
Using a projector $P_h$ with either 16 or 48 Fourier modes, we find 
an exponential rate convergence to begin immediately, with as few as 6 Chebyshev modes,
and no convergence with $\nc=4, 5$.   That $\nf\times\nc=16\times6$ modes
in $P_h$ would suffice at this higher Rayleigh number, but $\nf\times\nc=24\times6$ modes
in $P_h$ does not at $\Ra=2.5\times 10^7$ suggests an organizing effect of the flow as it becomes slightly more turbulent.

\section{Conclusions}  We have briefly recalled rigorous estimates on the global attractor of the Rayleigh-Benard system
in velocity form and interpreted how they influence the nudging algorithm for data assimilation studied here.  When those bounds
are inserted into the condition for synchronization they suggest that the resolution of the data might need to be exponentially small in 
$\Ra^{-1}$.  Even if one could reduce those bounds on the attractor to as small as $\OOO(1)$, the condition would still require a spatial resolution on the order of $\Ra^{-3/2}$.  

We then adapted the method for the vorticity formulation of the Rayleigh-Benard system, and presented numerical results to gauge 
the minimal amount of data (maximal resolution) needed for the algorithm, using a "zero-knowledge" initial condition, to synchronize at an exponential rate with a reference solution on the global attractor.  We looked at two Rayleigh number cases where the flow is nearly turbulent.  At $\Ra= 2.5\times 10^{7}$,
we found that data in as few as 6 Fourier $\times$ 8 Chebyshev modes or 24 Fourier $\times$ 7 Chebyshev modes is sufficient.  While at $\Ra= 5\times 10^{7}$, we find it is sufficient to take as few as 16 and 6 Fourier and Chebyshev modes respectively.

That synchronization is achievable at a resolution much more coarse than suggested by the rigorous analysis is consistent with experiments carried out at low Rayleigh number in \cite{Altaf} on the Rayleigh-Benard system in velocity form (which is also consistent with the computational study for the 2D NSE in \cite{Gesho_Olson_Titi}). A true comparison of our numerical results with analysis would require new estimates on the vorticity form of the system.
That should proceed using similar arguments to those in the velocity case, but would be considerably more complicated, so we leave it for future work.

\bigskip

\section*{Acknowledgements}

This work was initiated while the authors were visiting the Institute for Pure and Applied Mathematics (IPAM), which is supported by the National Science Foundation (NSF). The work of A.F. is supported in part by NSF grant DMS-1418911. The work of M.S.J. is supported in part by NSF grant DMS-1418911 and Leverhulme Trust grant VP1-2015-036. The work of E.S.T. is supported in part by the ONR grant N00014-15-1-2333.

\bigskip

\end{document}